\def\no{\if01}
\def\iftwelvept{\no}

\def\ifusepdf{\no}
\def\ifpsfont{\no}

\iftwelvept
\documentclass[leqno,12pt]{amsart}
\else
\documentclass[leqno,12pt]{amsart}
\fi
\usepackage{amssymb}
\usepackage{amscd}
\usepackage{latexsym}
\usepackage{verbatim}
\usepackage{mathrsfs}
\usepackage[all]{xy}

\ifusepdf
\usepackage{hyperref}
\else\fi
\ifpsfont
\usepackage[T1]{fontenc}
\usepackage{times}
\else\fi

\iftwelvept
\else\fi

\theoremstyle{plain}
\newtheorem{Theorem}{Theorem}[section]
\newtheorem*{TTheorem}{Theorem}

\newtheorem{Proposition}[Theorem]{Proposition}
\newtheorem{Lemma}[Theorem]{Lemma}
\newtheorem{Corollary}[Theorem]{Corollary}

\newtheorem{Claim}{Claim}[Theorem]

\theoremstyle{definition}

\newtheorem{Definition}[Theorem]{Definition}
\newtheorem{Remark}[Theorem]{Remark}
\newtheorem{Example}[Theorem]{Example}

\renewcommand{\theTheorem}{\arabic{section}.\arabic{Theorem}}

\newcommand{\AAAA}{{\mathbb{A}}}
\newcommand{\ZZ}{{\mathbb{Z}}}
\newcommand{\QQ}{{\mathbb{Q}}}
\newcommand{\RR}{{\mathbb{R}}}

\newcommand{\NN}{{\mathbb{N}}}
\newcommand{\GG}{{\mathbb{G}}}

\newcommand{\XXX}{{\mathscr{X}}}

\newcommand{\VVV}{{\mathscr{V}}}
\newcommand{\DDD}{{\mathscr{D}}}

\newcommand{\OO}{{\mathcal{O}}}
\newcommand{\LL}{{\mathcal{L}}}

\newcommand{\FF}{\mathcal{F}}

\newcommand{\Hom}{\operatorname{Hom}}

\newcommand{\Spec}{\operatorname{Spec}}

\newcommand{\mult}{\operatorname{mult}}

\newcommand{\Pic}{\operatorname{Pic}}

\newcommand{\Aut}{\operatorname{Aut}}

\newcommand{\Proof}{{\sl Proof.}\quad}
\newcommand{\PPP}{\mathcal{P}}
\newcommand{\QED}{{\unskip\nobreak\hfil\penalty50\quad\null\nobreak\hfil
{$\Box$}\parfillskip0pt\finalhyphendemerits0\par\medskip}}

\begin{document}
\title{Integral Chow rings of toric stacks}
\author{Isamu Iwanari}
\thanks{The author is supported by Grand-in-Aid for JSPS fellowships.}
\address{Department of Mathematics, Faculty of Science,
Kyoto University, Kyoto, 606-8502, Japan}
\email{iwanari@math.kyoto-u.ac.jp }
\begin{abstract}
The purpose of this paper is to prove that integral Chow rings of toric stacks
are naturally isomorphic to Stanley-Reisner rings.
\end{abstract}

\maketitle

\section*{Introductuon}
\renewcommand{\theTheorem}{\Alph{Theorem}}

Intersection theory with {\it integer coefficients} 
of Fulton-MacPherson type on smooth algebraic stacks
was developed by Kresch, Edidin-Graham and Totaro (\cite{EG},\cite{KR},\cite{TO}).
In particular, the integral Chow ring of a smooth stack
which has staratifications by quotient stacks (for example, quotient stacks)
was defined. The integral Chow ring of a smooth algebraic stack
is an interesting and
deep invariant of the stack which reflects
the geometric structure together with
the stacky structure on it.
It is challenging to compute them for interesting algebraic stacks.
Some examples of smooth Deligne-Mumford stacks were calculated.
For example, in \cite{EG} Edidin and Graham calculated
the integral Chow rings of the
moduli stacks of elliptic curves $\mathcal{M}_{1,1}$ and $\overline{\mathcal{M}}_{1,1}$.
Vistoli calculated the integral Chow ring of the moduli stack
of the curves of genus 2 (\cite{V2}).

The purpose of this paper is to compute the integral Chow rings
of toric stacks defined in (\cite{I}, \cite{I2}).
Our work of the presented paper
is motivated by the category-equivalence
between the 2-category of toric stacks and the
1-category of stacky fans (cf. \cite[Theorem1.2 and Theorem 1.4]{I2}).
Before stating  our main result, let us recall the result of Chow rings
of toric varieties.

\begin{TTheorem}[Fulton-Sturmfels, Danilov, Jurkiewicz]
Let $N=\ZZ^d$ be a lattice.
Let $\Delta$ be a non-singular fan in $N\otimes_{\ZZ}\RR$
and let $X_{\Delta}$ be the associated toric variety.
Let us denote by $A^*(X_{\Delta})$ the Chow ring of $X_{\Delta}$.
Then there exists a natural isomorphism of graded rings
\[
  \operatorname{(Stanley-Reisner\ ring\ of}
\Delta) \stackrel{\sim}{\to} A^*(X_{\Delta}).
\]
If  $\Delta$ is a simplicial fan and the base field is in characteristic zero,
then $A^*(X_{\Delta})\otimes_{\ZZ}\QQ$ has a ring structure and there exists
a natural isomorphism of graded rings
\[
 \operatorname{(Stanley-Reisner\ ring\ of}
\Delta)\otimes_{\ZZ}\QQ \stackrel{\sim}{\to} A^*(X_{\Delta})\otimes_{\ZZ}\QQ.
\]
\textup{(}See Definition~\ref{sr} for the definition of Stanley-Reisner rings
of $\Delta$.\textup{)}
\end{TTheorem}

Here we would like to invite the reader's
attention to the fact that  that if $\Delta$ is simplicial
and not non-singular, the (operational) Chow
group $A^k(X_{\Delta})$ for $k\ge 1$ could {\it differ
from the module of the ``degree $k$ part" of Stanley-Reisner ring
of $\Delta$} (cf. Example~\ref{exam}).
Furthermore in such case a somewhat surprising point is
that the Stanley-Reisner ring of $\Delta$
could be {\it nonzero
in degrees higher than the dimension of the toric variety $X_{\Delta}$}
(cf. Example~\ref{exam}). 
Since the Chow groups of an algebraic space in degrees higher than its
dimension are zero, thus Stanley-Reisner rings gives us a combinatorial phenomenon
which is unaccountable in the framework of schemes and algebraic spaces.

\vspace{1mm}

Now we state our main result.

\begin{TTheorem}[cf. Theorem~\ref{main2}]
Let $k$ be an algebraically closed field of characteristic zero.
Let $(\Delta,\Delta^0)$ be a stacky fan and
$\XXX_{(\Delta,\Delta^0)}$ the toric stack over $k$
associated to $(\Delta,\Delta^0)$.
Let $A^*(\XXX_{(\Delta,\Delta^0)})$ denote the
integral Chow ring of $\XXX_{(\Delta,\Delta^0)}$.
Suppose that rays in $\Delta$ span the vector space $N\otimes_{\ZZ}\RR$.
Then there exists a natural isomorphism of graded rings
\[
(\textup{Stanley-Reisner ring of}\ {(\Delta,\Delta^0)})\stackrel{\sim}{\to}A^*(\XXX_{(\Delta,\Delta^0)}).
\]
If $[\VVV(\sigma)]$ and $[\VVV(\tau)]$ are torus-invariant substacks in $\XXX_{(\Delta,\Delta^0)}$ which correspond to cones $\sigma$ and $\tau$ in $\Delta$ respectively
 \textup{(}cf. Section 2\textup{)},
we have
\begin{equation}
\tag*{(*)}
  [\VVV(\sigma)]\cdot [\VVV(\tau)]=[\VVV(\gamma)],
\end{equation}
in $A^*(\XXX_{(\Delta,\Delta^0)})$ when $\sigma$ and $\tau$ span $\gamma$.
\textup{(}See Definition~\ref{sr} for the definition of Stanley-Reisner rings
of $(\Delta,\Delta^0)$.\textup{)}
\end{TTheorem}

Since the Stanley-Reisner ring of a stacky fan
$(\Delta,\Delta_{\textup{can}}^0)$ with the canonical free-net $\Delta_{\textup{can}}^0$ coincides with the classical Stanley-Reisner ring of $\Delta$
(cf. Definition~\ref{sr}),
our result says that there exists an isomorphism
of graded rings between
the classical Stanley-Reisner ring of a simplicial fan $\Delta$
and $A^*(\XXX_{(\Delta,\Delta_{\textup{can}}^0)})$.

Here we explain how the usual relations on
 intersection product
of torus-invariant cycles on a simplicial toric variety $X_{\Sigma}$
(cf. \cite[page 100]{F}) 
are derived from that of the toric stack
$\XXX_{(\Delta,\Delta_{\textup{can}}^0)}$
in $A^*(\XXX_{(\Delta,\Delta_{\textup{can}}^0)})$.
There exists a coarse moduli map $\pi_{(\Delta,\Delta_{\textup{can}}^0)}:\XXX_{(\Delta,\Delta_{\textup{can}}^0)}\to X_{\Delta}$.
This functor defines the proper push-forward
$(\pi_{(\Delta,\Delta_{\textup{can}}^0)})_*:A^*(\XXX_{(\Delta,\Delta_{\textup{can}}^0)})\otimes_{\ZZ}\QQ \to A^*(X_{\Delta})\otimes_{\ZZ}\QQ$.
By \cite[Theorem 2.1.12 (ii)]{KR} and \cite[Proposition 6.1]{V},
$(\pi_{\Sigma})_*$ induces an isomorphism of groups.
Moreover, since a general stabilizer group of a toric stack is trivial,
$(\pi_{\Sigma})_*$ defines an isomorphism of rings (cf. \cite[(6.7)]{V}).
Thus the ring structure of $A^*(\XXX_{(\Delta,\Delta_{\textup{can}}^0)})\otimes_{\ZZ}\QQ$
yields that of $A^*(X_{\Delta})\otimes_{\ZZ}\QQ$.
The proper push-forward $(\pi_{(\Delta,\Delta_{\textup{can}}^0)})_*:
A^*(\XXX_{(\Delta,\Delta_{\textup{can}}^0)})\otimes_{\ZZ}\QQ {\to}
A^*(X_{\Delta})\otimes_{\ZZ}\QQ$
sends $[\VVV(\sigma)]$ to $\displaystyle \frac{1}{\mult (\sigma)}[V(\sigma)]$,
since the order of stabilizer group of a generic geometric point
on $\VVV(\sigma)$ is $\mult (\sigma)$
(cf. \cite[Proposition 4.13]{I}).  Here $\mult (\sigma)$ is the
multiplicity of $\sigma$, and
$V(\bullet)$ is the torus-invariant subvariety which corresponds to
a cone in $\Delta$.
Thus,
the relation (*)
in $A^*(\XXX_{(\Delta,\Delta_{\textup{can}}^0)})$ induces under $(\pi_{\Sigma})_*$ the relations
\[
  [V(\sigma)]\cdot [V(\tau)]= \displaystyle \frac{\mult (\sigma)\cdot \mult (\tau)}{\mult (\gamma)}[V(\gamma)]
\]
 in $A^*(X_{\Delta})\otimes_{\ZZ}\QQ$ if $\sigma$ and $\tau$ span $\gamma$.

\vspace{2mm}

There are more interesting points to notice.
It is known that the operational Chow group
of a complete toric variety is torsion-free (cf. \cite{FS}).
On the contrary, the integral Chow group of a complete toric stack
could have a lot of torsion elements.
Moreover, as noted above, it could be nonzero
in degrees higher than the dimension of the toric stack (cf. Example~\ref{exam}).
This is a substantial difference to the intersection theory on schemes
and algebraic spaces.

\vspace{1mm}

The presented result
can be also viewed as
an application of intersection theory with integral coefficients
(\cite{EG},\cite{KR},\cite{TO}) to toric stacks.
We hope that the reader finds our computation here shows a nice relation
between toric stacks and combinatorics.

\vspace{1mm}

The presented paper is organized as follows.
In section 1, for the computation of the integral Chow rings
(cf. \cite{EG}),
we obtain the quotient presentations of toric stacks defined
in \cite{I} and \cite{I2}. For this purpose,
we generalize the functor defined in \cite{C},
which represent a smooth toric variety (in characteristic zero),
to a certain groupoid.
In section 3, we present the proof of the main result.
Finally, we calculate some examples.

{\bf Notations And Conventions.}
Set $N=\ZZ^d$ and $M=\Hom_{\ZZ}(N,\ZZ)$.
Let $\langle\bullet,\bullet\rangle$ be the dual pairing.
Let $\Delta$ be a fan in $N_{\RR}=N\otimes_{\ZZ}\RR$
(we asumme that all fans are {\it finite} in this paper)
(cf. \cite{F}).
Denote by $\Delta(1)$ the set of rays.
Let us denote by $v_{\rho}$ the first lattice point on a ray
$\rho\in \Delta(1)$.
Finally, let $\Delta_{\textup{max}}$ denote the set of maximal
cones in $\Delta$.
A pair $(\Delta,\Delta^0)$ is called a {\it stacky fan}
if $\Delta$ is simplicial fan in $N_{\RR}$ and
$\Delta^0$ is a subset of $\Delta\cap N$
such that
for any cone $\sigma$ in $\Delta$, $\sigma\cap \Delta^0$ is
a sub-monoid of $\sigma\cap N$ which is
isomorphic to $\NN^r$ where $r=\dim \sigma$,
such that for any element $e \in \sigma\cap N$,
there exists a positive integer $n$
such that $n\cdot e \in \sigma\cap \Delta^0$.
The initial point of $\rho\cap \Delta^0$ is said to be
the {\it generator} of $\Delta^0$ on $\rho$.
Let $v_{\rho}$ denote the first lattice point of $\rho\cap N$ and
$n_{\rho}$ the initial point of $\rho\cap \Delta^0$.
The positive integer $l_{\rho}$
such that $l_{\rho}\cdot v_{\rho}=n_{\rho}$
is said to be the {\it level} of $\Delta^0$
on $\rho$. Notice that $\Sigma^0$ is completely determined by the levels
of $\Delta^0$ on rays of $\Delta$.
Each simplicial fan $\Delta$ has the canonical free-net
$\Delta^0_{\operatorname{can}}$, whose level on every ray in $\Delta$
is $1$.

Give a stacky fan $(\Delta,\Delta^0)$, we
have the associated toric stack $\XXX_{(\Delta,\Delta^0)}$
over a base scheme $S$.
If $S$ is the spectrum of a field $k$ of characteristic zero,
$\XXX_{(\Delta,\Delta^0)}$ is a smooth Deligne-Mumford stack that is
of finite type and separated over $k$.
For details, we refer to \cite[section 4]{I}, \cite{I2}.

\renewcommand{\theTheorem}{\arabic{section}.\arabic{Theorem}}
\renewcommand{\thesubsubsection}{\arabic{section}.\arabic{subsection}.\arabic{subsubsection}}

\section{$(\Delta,\Delta^0)$-collections}

In \cite{C}, given a fan $\Delta$ and a scheme $Y$,
Cox defines notions of {\it $\Delta$-collections} on $Y$ and equivalences
between them.
Then he showed that the functor$
F_{\Delta}:Y \mapsto \{ \Delta\textup{-collections\ on } Y\}/\sim$
represents the toric variety $X_{\Delta}$ if $\Delta$ is {\it non-singular}.
If $\Delta$ is singular, unfortunately, the functor of $\Delta$-collections
fails to represent the toric variety $X_{\Delta}$.
The aim of this section is to generalize the notion of $\Delta$-collections
and \cite[Theorem 1.1]{C}
and give a quotient presentation for a toric stack $\XXX_{(\Delta,\Delta^0)}$
defined in \cite{I} (cf. Corollary\ref{quot}).

\vspace{1mm}

\begin{Definition}
A {\it $(\Delta,\Delta^0)$-collection} on $Y$ is a triple
\[
(\{\LL_{\rho}\}_{\rho\in \Delta(1)},\{u_{\rho}\}_{\rho\in \Delta(1)}, 
\{ c_m:\otimes_{\rho}\LL_{\rho}^{\otimes\langle m,n_{\rho}\rangle}\stackrel{\sim}{\to} \OO_{Y}\}_{m\in M}),
\]
where
$\LL_{\rho}$ is an invertible sheaf on $Y$,
$u_{\rho}\in H^0(Y,\LL_{\rho})$ and $c_m$ is an isomorphism of invertible sheaves with the following additional properties:
\begin{enumerate}
\item $c_m\otimes c_{m'}=c_{m+m'}$ for \textup{all} $m, m'\in M$.
\item The homomorphism $\Sigma_{\sigma\in \Delta_{\textup{max}}}\otimes_{\rho \notin \sigma(1)}u_{\rho}^*:\bigoplus_{\sigma\in \Delta_{\textup{max}}}\otimes_{\rho\notin \sigma(1)}\LL_{\rho}^{-1}\to \OO_Y$ is surjective
where $u_{\rho}^*:\LL_{\rho}^{-1}\to \OO_Y$ is induced by
the homomorphism $u_{\rho}:\OO_Y\to \LL_{\rho}, (1\mapsto u_{\rho})$.
\end{enumerate}
If $f:Y'\to Y$ is a morphism of schemes,
then we define the pullback 
\[
f^*(\{\LL_{\rho}\}_{\rho\in \Delta(1)},\{u_{\rho}\}_{\rho\in \Delta(1)}, 
\{ c_m:\otimes_{\rho}\LL_{\rho}^{\otimes\langle m,n_{\rho}\rangle}\stackrel{\sim}{\to} \OO_{Y}\}_{m\in M}),
\]
to be $(\{f^*\LL_{\rho}\}_{\rho\in \Delta(1)},\{f^*u_{\rho}\}_{\rho\in \Delta(1)}, 
\{ f^*c_m:\otimes_{\rho}f^*\LL_{\rho}^{\otimes\langle m,n_{\rho}\rangle}\stackrel{\sim}{\to} \OO_{Y'}\}_{m\in M})$.
\end{Definition}

\begin{Remark}
Let $\Delta$ be a non-singular fan and $\Delta_{\textup{can}}^0$
the canonical free-net.
Then $\Delta$-collections in \cite{C} are $(\Delta,\Delta_{\textup{can}}^0)$-collections,
and vice versa.

If $\Delta$ is the empty set, a $(\Delta,\Delta^0)$-collection on $Y$ is
just a collection $\{c_{m}:\OO_{Y}\stackrel{\sim}{\to}\OO_{Y}\}_{m\in M}$
such that $c_{m+m'}=c_{m}\otimes c_{m'}$ for all $m, m'\in M$.
This canonically corresponds to a $Y$-valued points of $\Spec \ZZ[M]$.
\end{Remark}

\begin{Definition}
\label{mor}
If $\Delta$ is a non-empty set, a morphism 
\[
(\{\LL_{\rho}\}_{\rho\in \Delta(1)},\{u_{\rho}\}_{\rho\in \Delta(1)}, 
\{ c_m\}_{m\in M})\to
(\{\LL_{\rho}'\}_{\rho\in \Delta(1)},\{u_{\rho}'\}_{\rho\in \Delta(1)}, 
\{ c_m'\}_{m\in M})
\]
between $(\Delta,\Delta^0)$-collections on $Y$ is
a set of isomorphisms $\{ \phi_{\rho}:\LL_{\rho}\stackrel{\sim}{\to}
\LL_{\rho}'\}_{\rho\in \Delta(1)}$ indexed by $\Delta(1)$,
such that:
\begin{enumerate}
\item $\phi_{\rho}(u_{\rho})=u_{\rho}'$
\item For every $m\in M$ the diagram
\[
\xymatrix{
\otimes_{\rho}\LL_{\rho}^{\otimes\langle m,n_{\rho}\rangle} \ar[r]^(0.6){c_m} \ar[d]_{\otimes_{\rho}\phi_{\rho}^{\otimes\langle m,n_{\rho}\rangle}} & \OO_Y \\
\otimes_{\rho}\LL_{\rho}^{'\otimes\langle m,n_{\rho}\rangle} \ar[ur]_{c_m'}\\
}
\]
commutes.

When $\Delta$ is the empty set,
the set of morphisms from $\{ c_m:\OO_Y\stackrel{\sim}{\to}\OO_Y\}$
to $\{ c_m':\OO_Y\stackrel{\sim}{\to}\OO_Y\}$
is $\{ \textup{Idetity} \}$ if $c_m=c_{m}'$ for all $m\in M$
and is the empty set if otherwise.
\end{enumerate}
\end{Definition}

Let $S$ be a scheme, and define a fibered category
\[
\FF_{(\Delta,\Delta^0)}\to (S\textup{-schemes})
\]
as follows.
The objects of $\FF_{(\Delta,\Delta^0)}$ over a $S$-scheme $Y$
are $(\Delta,\Delta^0)$-collections on $Y$.
A morphism between two objects in $\textup{Ob}(\FF_{(\Delta,\Delta^0)})(Y)$
is a morphism of $(\Delta,\Delta^0)$-collections on $Y$.
With the natural notion of pullbacks, $\FF_{(\Delta,\Delta^0)}$
is a fibered category over $(S\textup{-schemes})$.
By fppf descent theory for quasi-coherent sheaves,
$\FF_{(\Delta,\Delta^0)}$ is a stack with respect to fppf topology.

\begin{Theorem}
\label{m}
Let $S$ be the spectrum of an algebraically closed field $k$ of characteristic
zero. Let $\XXX_{(\Delta,\Delta^0)}$ be the toric stack (over $k$) associated to $(\Delta,\Delta^0)$ \textup{(cf. \cite{I2})}.
Then there exists an isomorphism of stacks
\[
\FF_{(\Delta,\Delta^0)}\stackrel{\sim}{\longrightarrow} \XXX_{(\Delta,\Delta^0)}.
\]
\end{Theorem}

The proof of Theorem~\ref{m} proceeds in several steps.

\begin{Lemma}
\label{sub}
To prove Theorem~\ref{main} it suffices to show the followings:
\begin{enumerate}
\renewcommand{\labelenumi}{(\alph{enumi})}

\item The stack $\FF_{(\Delta,\Delta^0)}$ is a smooth Deligne-Mumford stack
of finite type and separated over $k$.
\item Let $\pi:\FF_{(\Delta,\Delta^0)}\to F$ be a coarse moduli space
for $\FF_{(\Delta,\Delta^0)}$ \textup{(}see for example, \cite[Definition 2.4]{hom}\textup{)}. Then $F$ is isomorphic to the toric variety
$X_{\Delta}$, and $\pi^{-1}(\Spec k[M])\to \Spec k[M]\subset \XXX_{(\Delta,\Delta^0)}$ is an isomorphism. If $D_{\rho}$ is the toric divisor of $X_{\Delta}$
corresponding to $\rho$, the order of the stabilizer group of a geometric point on the
generic point of $\pi^{-1}(D_{\rho})_{\textup{red}}$ is the level of
$\Delta^0$ on $\rho$.

\end{enumerate}
\end{Lemma}

\Proof
By proof of \cite[Theorem 1.3]{I2} and above two assumptions,
$\FF_{(\Delta,\Delta^0)}$ is isomorphic to
the toric stacks $\XXX_{(\Sigma,\Sigma^0)}$ for some stacky fan $(\Sigma,\Sigma^0)$ (Note that if the coarse moduli space is a
toric variety the assumption of an existence of a torus action
in \cite[Theorem 1.3]{I2} is not necessary for our purpose).
We have to show $(\Delta,\Delta^0)=(\Sigma,\Sigma^0)$.
Since $F=X_{\Delta}$ we have $\Delta=\Sigma$.
Let $\DDD_{\rho}\subset \XXX_{(\Delta,\Delta^0)}$
(resp. $\DDD_{\rho}'\subset \XXX_{(\Delta,\Sigma^0)}$)
the toric divisor corresponding to $\rho$ in $\XXX_{(\Delta,\Delta^0)}$
(resp. $\XXX_{(\Delta,\Sigma^0)}$).
Then notice that $\XXX_{(\Delta,\Delta^0)}$
is isomorphic to $\XXX_{(\Delta,\Sigma^0)}$
if and only if for any ray $\rho\in\Delta$
both stabilizer groups of geometric points on the generic points
of $\DDD_{\rho}$ and $\DDD_{\rho}'$ have the same orders.
Since the order of the stabilizer group of a geometric point on the generic points
of $\DDD_{\rho}$ equals to the level of $\Delta^0$ on $\rho$,
we have $\Delta^0=\Sigma^0$.
\QED

\subsection{A smooth Deligne-Mumford stack of finite type and separated over $k$.}

\vspace{1mm}
In this subsection, we prove (a) in Lemma~\ref{sub}.
Unless stated otherwise, we work over $k$.
Consider a collection
\[
 P=( \{\alpha_{\rho}:\OO_{Y,\rho}\stackrel{\sim}{\to}\OO_Y\}_{\rho\in \Delta(1)}, \{ u_{\rho}\in H^0(Y,\OO_{Y,\rho}) \}_{\rho\in \Delta(1)}, \{ c_r:\OO_Y\stackrel{\sim}{\to}\OO_Y\}_{r\in R})
\]
where $\alpha_{\rho}:\OO_{Y,\rho}\stackrel{\sim}{\to}\OO_Y$
is an isomorphism of invertible sheaves, such that $\Sigma_{\sigma\in \Delta_{\textup{max}}}\otimes_{\rho \notin \sigma(1)}u_{\rho}^*:\bigoplus_{\sigma\in \Delta_{\textup{max}}}\otimes_{\rho\notin \sigma(1)}\OO_{Y,\rho}^{-1}\to \OO_Y$ is surjective,
and $c_{r+r'}=c_{r}\otimes c_{r'}$ for all $r,r'\in R$.
Here $R=\Hom_{\textup{groups}}(N/(\Sigma_{\rho}\ZZ\cdot v_{\rho}),\ZZ)\subset M$ where
$v_{\rho}$ is the first lattice point of $\rho$.
We shall refer to such collections as {\it linear $\Delta$-collections}.

Let $P_i=( \{\alpha_{\rho}^{(i)}:\OO_{Y,\rho}^{(i)}\stackrel{\sim}{\to} \OO_Y\}_{\rho\in \Delta(1)}, \{ u_{\rho}^{(i)}\in H^0(Y,\OO_{Y,\rho}) \}_{\rho\in \Delta(1)}, \{ c_r^{(i)}:\OO_Y\stackrel{\sim}{\to}\OO_Y\}_{r\in R})$ (for $i=1,2$) be linear $\Delta$-collections.
We say that $P_1$ is equivalent to $P_2$
if there exists an isomorphism $\phi_{\rho}:\OO_{Y,\rho}^{(1)}\to \OO_{Y,\rho}^{(2)}$
such that $\phi_{\rho}(u_{\rho}^{(1)})=u_{\rho}^{(2)}$ and $\alpha_{\rho}^{(2)}\circ \phi_{\rho}=\alpha_{\rho}^{(1)}$ for all $\rho\in \Delta(1)$
and $c_{r}^{(1)}=c_{r}^{(2)}$ for all $r\in R$.
%If $\Delta$ is a non-empty set and $c_{r}^{(1)}=c_{r}^{(2)}$ for all $r\in R$,
%a morphism $P\to P'$ is a collection of isomorphisms of invertible sheaves
%$\{ \phi_{\rho}:\OO_{Y,\rho}\stackrel{\sim}{\to}\OO_{Y,\rho}'\}_{\rho\in\Delta(%1)}$ such that $\phi_{\rho}(u_{\rho})=u_{\rho}'$ for all $\rho\in\Delta(1)$.
%If $\Delta$ is the empty set and $c_{r}^{(1)}=c_{r}^{(2)}$ for all $r\in R$, a %morphism $P\to P'$ is an ``identity''.
%If otherwise, the set of morphisms is the empty set.
The following functor
\[
L_{\Delta}:(k\textup{-schemes})\to (\textup{Sets}),\ \ Y\mapsto \{ \textup{linear}\ \Delta\textup{-collections}\}/\sim
\]
(here $\sim$ denotes the equivalence relation) is representable by
a quasi-affine scheme
$\bigcup_{\sigma\in \Delta_{\textup{max}}}\Spec k[x_{\rho},\rho\in \Delta(1)]_{\Pi_{\rho\notin \sigma(1)}x_{\rho}}\otimes k[R]\subset \Spec k[x_{\rho},\rho\in \Delta(1)]\otimes k[R]=\AAAA^{\Delta(1)}\times \Spec k[R]$.
Define a group scheme
\begin{align*}
G_{(\Delta,\Delta^0)} &:=\{ (a_{\rho} )_{\rho\in\Delta(1)}\in \GG_m^{\Delta(1)}|\Pi_{\rho}a_{\rho}^{\langle m,n_{\rho}\rangle}=1\ \textup{for any}\ m\in M\} \\
&=\Spec k[x_{\rho}^{\pm 1},\rho\in\Delta(1)]/(\Pi_{\rho}x_{\rho}^{\langle m,n_{\rho}\rangle}-1)_{m\in M}.
\end{align*}
There exists a natural action $a:L_{\Delta}\times G_{(\Delta,\Delta^0)}\to L_{\Delta}$ of $G_{(\Delta,\Delta^0)}$ on $L_{\Delta}$
defined by
\[
\{\{\alpha_{\rho}:\OO_{Y,\rho}\stackrel{\sim}{\to} \OO_Y\},\{u_{\rho} \},\{ c_r:\OO_Y\stackrel{\sim}{\to}\OO_Y\}\} \cdot (a_{\rho}):=\{\{\alpha_{\rho}:\OO_{Y,\rho}\stackrel{\sim}{\to}\OO_Y\},\{u_{\rho}\cdot  a_{\rho}\}, \{ c_r:\OO_Y\stackrel{\sim}{\to}\OO_Y\}\}.
\]
Since $G_{(\Delta,\Delta^0)}$ is a smooth group scheme (characteristic zero),
the quotient stack $[L_{\Delta}/G_{(\Delta,\Delta^0)}]$ is a smooth algebraic stack
of finite type over $k$ (cf. \cite[proposition 10.13.1]{LM}). 

Next we show that $[L_{\Delta}/G_{(\Delta,\Delta^0)}]$ has a finite diagonal.
To this aim, clearly we may suppose that rays in $\Delta$
span the vector space $N\otimes_{\ZZ}\RR$.
For a cone $\sigma\in \Delta$ set $U_{\sigma}=\Spec k[x_{\rho},\rho\in \Delta(1)]_{\Pi_{\rho\notin \sigma(1)}x_{\rho}}\subset L_{\Delta}=\cup_{\delta\in \Delta_{\textup{max}}}U_{\delta}$.
The restriction $U_{\sigma}\times G_{(\Delta,\Delta^0)}\to U_{\sigma}\times U_{\sigma}$
of the morphism $L_{\Delta}\times G_{(\Delta,\Delta^0)}\stackrel{\textup{pr}_{1}\times a}{\longrightarrow} L_{\Delta}\times L_{\Delta}$ is defined by
\[
\begin{CD}
k[x_{\rho},\rho\in \Delta(1)]_{\Pi_{\rho\notin \sigma(1)}x_{\rho}}\otimes
k[x_{\rho},\rho\in \Delta(1)]_{\Pi_{\rho\notin \sigma(1)}x_{\rho}} \\
@VVV \\
k[x_{\rho},\rho\in \Delta(1)]_{\Pi_{\rho\notin \sigma(1)}x_{\rho}}\otimes k[x_{\rho}^{\pm 1},\rho\in\Delta(1)]/(\Pi_{\rho}x_{\rho}^{\langle m,n_{\rho}\rangle}-1)_{m\in M}, \\
\end{CD}
\]
\[
\scriptsize{x_{\rho}\otimes 1\mapsto x_{\rho}\otimes 1 ,\ 1\otimes x_{\rho}\mapsto x_{\rho}\otimes x_{\rho}}.
\]
It suffices to show that $k[x_{\rho},\rho\in \Delta(1)]_{\Pi_{\rho\notin \sigma(1)}x_{\rho}}\otimes k[x_{\rho}^{\pm 1},\rho\in\Delta(1)]/(\Pi_{\rho}x_{\rho}^{\langle m,n_{\rho}\rangle}-1)_{m\in M}$ is a finite $k[x_{\rho},\rho\in \Delta(1)]_{\Pi_{\rho\notin \sigma(1)}x_{\rho}}\otimes
k[x_{\rho},\rho\in \Delta(1)]_{\Pi_{\rho\notin \sigma(1)}x_{\rho}}$ algebra.
By an easy observation, it is reduced to show that $k[x_{\rho}^{\pm 1}, \rho\in \sigma(1)]/(\Pi_{\rho\in \sigma(1)}x_{\rho}^{\langle m,n_{\rho}\rangle}-1)_{m\in M}$ is a finite $k$-algebra. This follows from the fact that
$\sigma$ is simplicial. Thus $[L_{\Delta}/G_{(\Delta,\Delta^0)}]$ has a finite diagonal,
in particular separated over $k$. Moreover in characteristic zero,
it is a Deligne-Mumford stack.

Now we construct an isomorphism from $[L_{\Delta}/G_{(\Delta,\Delta^0)}]$ to $\FF_{(\Delta,\Delta^0)}$. This part is non-canonical.
Choose a splitting $M\cong R\oplus R'$ where $R'$ is a finitely generated
free abelian group.
Define a morphism $L_{\Delta}\to \FF_{(\Delta,\Delta^0)}$
by $( \{\alpha_{\rho}:\OO_{Y,\rho}\stackrel{\sim}{\to}\OO_Y\}_{\rho\in \Delta(1)}, \{ u_{\rho} \}_{\rho\in \Delta(1)}, \{ c_r:\OO_Y\stackrel{\sim}{\to}\OO_Y\}_{r\in R})\mapsto (\{\OO_{Y,\rho},u_{\rho}\}_{\rho\in\Delta(1)},\{ c_{\operatorname{pr}_1(m)}\otimes\otimes_{\rho}\alpha_{\rho}^{\otimes\langle m,n_{\rho}\rangle}:\OO_Y\otimes_{\rho}\OO_{Y,\rho}^{\otimes\langle m,n_{\rho}\rangle}\to\OO_Y \}_{m\in M})$,
where $\operatorname{pr}_1:M=R\oplus R'\to R$ is the first projection.
Consider the diagram 
$L_{\Delta}\times G_{(\Delta,\Delta^0)}\rightrightarrows L_{\Delta} \to \FF_{(\Delta,\Delta^0)}$.
Then we have a natural 2-isomorphism between the both two composites.
Thus there exists a morphism $z:[L_{\Delta}/G_{(\Delta,\Delta^0)}]\to \FF_{(\Delta,\Delta^0)}$.
We can easily observe that every $(\Delta,\Delta^0)$-collection on a $k$-scheme $Y$
has, fppf locally on $Y$, the form
$(\{ \OO_{Y},u_{\rho}\}_{\rho\in \Delta(1)},\{ c_{m}:\otimes_{\rho}\OO_{Y}^{\otimes\langle m,n_{\rho}\rangle}\to\OO_Y \}_{m\in M})$
such that $c_{r'}$ is an identity for all $r'\in R'$.
Since $[L_{\Delta}/G_{(\Delta,\Delta^0)}]$ and $\FF_{(\Delta,\Delta^0)}$
are fppf stacks (cf. \cite[(10.7)]{LM}), the functor $z$ is essentially surjective.

To prove the fully faithfulness,
note first that we may work fppf locally on $Y$
and put two linear $\Delta$-collections
$C_1:=( \{\alpha_{\rho}:\OO_{Y,\rho}\stackrel{\sim}{\to}\OO_Y\}_{\rho\in \Delta(1)}, \{ u_{\rho} \}_{\rho\in \Delta(1)}, \{ c_r:\OO_Y\stackrel{\sim}{\to}\OO_Y\}_{r\in R})$
and
$C_2:=( \{\alpha_{\rho}:\OO_{Y,\rho}\stackrel{\sim}{\to}\OO_Y\}_{\rho\in \Delta(1)}, \{ u_{\rho}' \}_{\rho\in \Delta(1)}, \{ c_r':\OO_Y\stackrel{\sim}{\to}\OO_Y\}_{r\in R})$.
If $c_r=c_r'$ for all $r\in R$, the set of morphisms from $C_1$ to $C_2$ in $[L_{\Delta}/G](Y)$ is
\[
\{(\lambda_{\rho})_{\rho}\in G(Y)|\ \lambda_{\rho}\cdot u_{\rho}=u_{\rho}'\ \textup{for any}\ \rho\in \Delta(1) \}.
\]
In this case, the set of homomorphism from $z(C_1)=(\{ \OO_{Y,\rho},u_{\rho}\}_{\rho\in \Delta(1)},\{ c_{\operatorname{pr}_1(m)}\otimes\otimes_{\rho}\alpha_{\rho}^{\otimes\langle m,n_{\rho}\rangle}:\OO_Y\otimes_{\rho}\OO_{Y,\rho}^{\otimes\langle m,n_{\rho}\rangle}\to\OO_Y \}_{m\in M})$ to $z(C_2)=(\{ \OO_{Y,\rho},u_{\rho}'\}_{\rho\in \Delta(1)},\{ c_{\operatorname{pr}_1(m)}'\otimes\otimes_{\rho}\alpha_{\rho}^{\otimes\langle m,n_{\rho}\rangle}:\OO_Y\otimes_{\rho}\OO_{Y,\rho}^{\otimes\langle m,n_{\rho}\rangle}\to\OO_Y \}_{m\in M})$ in $\FF_{(\Delta,\Delta^0)}(Y)$
can be canonically identified
with the set
\[
\{ (\lambda_{\rho})_{\rho}\in \GG_m^{\Delta(1)}(Y)| \lambda_{\rho}u_{\rho}=u_{\rho}'\ \textup{for any}\ \rho\in \Delta(1), \Pi_{\rho}\lambda_{\rho}^{\langle m,n_{\rho}\rangle}=1\ \textup{for any}\ m\in M \}
\]
by $\{ \lambda_{\rho}:\OO_{Y,\rho}\to \OO_{Y,\rho}, 1\mapsto \lambda_{\rho}\}_{\rho}$, and bijectively corresponds to the set of morphisms from
$C_1$ to $C_2$ in $[L_{\Delta}/G](Y)$.
If otherwise, both $\Hom_{[L_{\Delta}/G](Y)}(C_1,C_2)$
and $\Hom_{\FF_{(\Delta,\Delta^0)}(Y)}(z(C_1),z(C_2))$
are the empty sets. Thus $z$ is fully faithful.
Therefore $\FF_{(\Delta,\Delta^0)}$ is a smooth Deligne-Mumford stack
of finite type and separated over $k$.

\begin{Remark}
The above argument also implies that
the stack $\FF_{(\Delta,\Delta^0)}$ over a general scheme
is algebraic. Namely, we conclude that:
\end{Remark}

\begin{Proposition}
Let $\FF_{(\Delta,\Delta^0)/\ZZ}$ denote the stack of $(\Delta,\Delta^0)$-collections
over $\Spec \ZZ$.
Then $\FF_{(\Delta,\Delta^0)/\ZZ}$ is an Artin stack
of finite type over $\ZZ$ with finite diagonal.
Let $r$ be the order of the torsion subgroup of the cokernel
of $M\to \oplus_{\rho\in \Delta(1)}\ZZ\cdot e_{\rho}$ defined by $m\mapsto \Sigma_{\rho}\langle m,n_{\rho}\rangle\cdot e_{\rho}$.
If $r$ is invertible on a scheme $S$, $\FF_{(\Delta,\Delta^0)/\ZZ}\times_{\ZZ}S$ is smooth over $S$.
\end{Proposition}

\begin{Remark}
Let $(\Delta,\Delta^0_{\textup{can}})$ be a stacky fan with the canonical free-net such that $\Delta$ is non-singular.
Let $X_{\Delta}$ (resp. $\XXX_{(\Delta,\Delta^0_{\textup{can}})}$) be the toric variety
(resp. the toric stack)
over $\ZZ$ associated to $\Delta$ (resp. $(\Delta,\Delta^0_{\textup{can}})$)
(the definition of toric stacks (\cite{I}, \cite{I2}) works over arbitrary base schemes).
While $\XXX_{(\Delta,\Delta^0_{\textup{can}})}$ is isomorphic to the toric variety $X_{\Delta}$
over $\ZZ$,  it is not clear whether or not
$\FF_{(\Delta,\Delta^0)/\ZZ}$ is isomorphic to $X_{\Delta}$ over $\ZZ$.

\end{Remark}

\subsection{The coarse moduli space for $\FF_{(\Delta,\Delta^0)}$}
In this subsection, we prove (b) in Lemma~\ref{sub}.
Clearly, we may suppose that rays in $\Delta$ span
the vector space $N\otimes_{\ZZ}\RR$.
Set $G:=G_{(\Delta,\Delta^0)}$.

First by imitating the proof of \cite[Theorem 2.1]{C2},
we see that geometric quotient (in the sense of Mumford \cite{M})
of $L_{\Delta}\times G\stackrel{a}{\to} L_{\Delta}$ is
a toric variety $X_{\Delta}$ associated to $\Delta$.
We will show that the toric variety $X_{\Delta}$
is a coarse moduli space for $[L_{\Delta}/G]\cong \FF_{(\Delta,\Delta^0)}$.
To this aim, by \cite[Theorem 2.6 (iii)]{Hom},
it suffices only to prove that $q:[L_{\Delta}/G]\to X_{\Delta}$ is proper
with the property $\OO_{X_{\Delta}}\cong q_*\OO_{[L_{\Delta}/G]}$
and gives a bijection on geometric points.
Note first that $X_{\Delta}$ is a geometric quotient and thus
$q$ induces a bijection on geometric points.
Moreover $\OO_{X_{\Delta}}\cong q_*\OO_{[L_{\Delta}/G]}$.
The properness of $q$ follows from the fact that $q$ is a universal submersion,
in particular universal closed map
(It is easy to see that $q$ is separated and of finite type).
Therefore the coarse moduli space for $\FF_{(\Delta,\Delta^0)}$
is a toric variety $X_{\Delta}$.

The quasi-affine scheme $L_{\Delta}$ contains
an algebraic torus $\Spec k[x_{\rho}^{\pm 1},\rho\in \Delta(1)]$
as a dense open set.
The restriction $\Spec k[x_{\rho}^{\pm 1},\rho\in \Delta(1)]\times G\to \Spec k[x_{\rho}^{\pm 1},\rho\in \Delta(1)]$
of the action $a:L_{\Delta}\times G\to L_{\Delta}$
is free and the quotient stack $[\Spec k[x_{\rho}^{\pm 1},\rho\in \Delta(1)]/G]$
is isomorphic to $\Spec k[M]$ by \cite[Claim 4.1.1]{I}.
The torus $\Spec k[M]\subset \FF_{(\Delta,\Delta^0)}$
is identified with the torus in $X_{\Delta}$
via the coarse moduli map.

Next fix a ray $\rho_{0}$. We will show
that if $D_{\rho_0}$ denotes the toric divisor of $X_{\Delta}$
corresponding to $\rho_0$, the order of the stabilizer group of a geometric point on the
generic point of $\pi^{-1}(D_{\rho})_{\textup{red}}$ is the level of
$\Delta^0$ on $\rho_0$.
Let $K$ be an algebraic closure of the function field of
$\pi^{-1}(D_{\rho})_{\textup{red}}$ and $p:\Spec K\to \FF_{(\Delta,\Delta^0)}$
a geometric point over the
generic point of $\pi^{-1}(D_{\rho})_{\textup{red}}$.
Taking account of the smooth atlas $L_{\Delta}\to \FF_{(\Delta,\Delta^0)}$,
the $(\Delta,\Delta^0)$-collection on $\Spec K$ determined by
$p:\Spec K\to \FF_{(\Delta,\Delta^0)}$ has the form
$(\{ \OO_{\Spec K,\rho},u_{\rho}\}_{\rho}, \{ c_m:\otimes_{\rho}\OO_{\Spec K,\rho}^{\langle m,n_{\rho}\rangle}\to \OO_{\Spec K}\}_{m\in M})$
such that $u_{\rho_o}=0$ and $u_{\rho}\neq 0$ if $\rho\neq\rho_0$.
It remains to calculate the order of the group of automorphisms
of this $(\Delta,\Delta^0)$-collection.
To this end, note that we may assume
$\rho_0=\RR_{\ge 0}\cdot(1,0,\ldots,0)\subset \ZZ^d=N$.
Furthermore if $\{ \phi_{\rho}\}_{\rho}$ is an automorphism,
$\phi_{\rho}$ ($\rho\neq\rho_0$) should be an identity 1 since $u_{\rho}\neq0$.
Let $l$ be the level of $\Delta^0$ on $\rho_0$.
Then $n_{\rho_0}=(l,0,\ldots,0)$.
By the above argument,
$\{\phi_{\rho_0}:\OO_{\Spec K,\rho_0}\to \OO_{\Spec K,\rho_0},\phi_{\rho}=1\ (\rho\neq \rho_0)\}$ gives an automorphism exactly when
for any $m\in M$, the diagram
\[
\xymatrix{
\OO_{\Spec K,\rho_0}^{\otimes\langle m,n_{\rho_0}\rangle} \ar[r]^(0.6){c_m} \ar[d]_{\phi_{\rho_0}^{\otimes\langle m,n_{\rho_0}\rangle}} & \OO_{\Spec K} \\
\OO_{\Spec K,\rho_0}^{\otimes\langle m,n_{\rho_0}\rangle} \ar[ur]_{c_m}\\
}
\]
commutes.
Thus $\phi_{\rho_0}:\OO_{\Spec K,\rho_0}\to \OO_{\Spec K,\rho}$
sends $1\in\OO_{\Spec K,\rho_0}$ to $\zeta\in \OO_{\Spec K,\rho_0}$
for some $\zeta\in \mu_{l}(K)$, and the automorphism group
is isomorphic to $\mu_l(K)\cong \ZZ/l\ZZ$.

\vspace{1mm}

We complete the proof of Theorem~\ref{m}.\QED

\begin{Corollary}
\label{quot}
\label{main}
Let $k$ be an algebraically closed field $k$ of characteristic
zero. Let $\XXX_{(\Delta,\Delta^0)}$ be the toric stack (over $k$) associated to $(\Delta,\Delta^0)$ \textup{(cf. \cite{I2})}.
Then there exists an isomorphism of stacks over $k$
$[L_{\Delta}/G_{(\Delta,\Delta^0)}]\stackrel{\sim}{\longrightarrow}\FF_{(\Delta,\Delta^0)}\stackrel{\sim}{\longrightarrow} \XXX_{(\Delta,\Delta^0)}$.
\end{Corollary}

\section{Integral Chow Rings of Toric Stacks}

In this section, we calculate integral Chow rings
of toric stacks.
We shall use notation similar to section 1, and from now on we assume that the base field $k$ is {\it an algebraically closed field of characteristic zero}.
Our computation is based on
intersection theory on stacks due to
 Kresch, Edidin-Graham, and Totaro.
For details, we refer to \cite{EG} \cite{KR} \cite{TO}.
For a Deligne-Mumford stack $\mathcal{X}$, we denote
by $A_{i}(\mathcal{X})$ the $i$-th {\it integral Chow group} of $\mathcal{X}$
(cf. \cite[section 5.3]{EG}, \cite[section 2.1]{KR}).

Let us fix some notations.
If $\Delta$ is a fan (resp. a stacky fan),
then for a cone $\delta\in \Delta$ we denote
by $V(\delta)$ (resp. $\VVV(\delta)$)
the torus-invarint cycle on the toric variety $X_{\Delta}$
(resp. the toric stack $\XXX_{(\Delta,\Delta^0)}$),
which corresponds to $\delta$.
If $\pi_{(\Delta,\Delta^0)}:\XXX_{(\Delta,\Delta^0)}\to X_{\Delta}$
denotes a coarse moduli map, then the cycle $\VVV(\delta)$
defines to be the reduced cycle $\pi_{(\Delta,\Delta^0)}^{-1}(V(\delta))_{\textup{red}}$.
For a ray $\rho\in \Delta$,
if no confusion seems likely to arise,
we may write $D_{\rho}$ (resp. $\DDD_{\rho}$)
for the torus-invarinat divisor $V(\rho)$ (resp. $\VVV(\rho)$).

\vspace{1mm}

Given a stacky fan $(\Delta,\Delta^0)$,
let us define {\it the Stanley-Reisner ring of $(\Delta,\Delta^0)$}.

\begin{Definition}[Stanley-Reisner ring]
\label{sr}
Let $(\Delta,\Delta^0)$ be a stacky fan.
Consider a polynomial ring $\ZZ[\DDD_{\rho}, \rho\in \Delta(1)]$.
Let $I_{(\Delta,\Delta^0)}$ be the ideal generated by
the linear forms $\Sigma_{\rho\in \Delta(1)}\langle m,n_{\rho} \rangle \DDD_{\rho}$, as $m$ ranges over $M$.
Here $n_{\rho}$ is the generator of $\Delta^0$ on $\rho$.
Let $J_{\Delta}$ be the ideal generated by the monomials
$\DDD_{\rho_1}\cdots \DDD_{\rho_s}$
such that $\langle \rho_{1},\ldots, \rho_{s}\rangle \notin \Delta$.
We define {\it the Stanley-Reisner ring $\textup{SR}_{(\Delta,\Delta^0)}$
of $(\Delta,\Delta^0)$} to be
\[
\ZZ[\DDD_{\rho}, \rho\in \Delta(1)]/(I_{(\Delta,\Delta^0)}+J_{\Delta}).
\]
For a simplicial fan $\Delta$ let us denote by $\Delta_{\textup{can}}^0$ the canonical free-net
and set $I_{\Delta}:=I_{(\Delta,\Delta_{\textup{can}}^0)}$.
Then the classical Stanley-Reisner ring $\textup{SR}_{\Delta}$ of $\Delta$
is defined to be
\[
\ZZ[\DDD_{\rho}, \rho\in \Delta(1)]/(I_{\Delta}+J_{\Delta}).
\]
For the notion of Stanley-Reisner rings of fans and polytopes, see (\cite[Section 5.2]{F},
\cite[10.7]{D}).
\end{Definition}

\begin{Theorem}
\label{main2}
Let $(\Delta,\Delta^0)$ be a stacky fan and
$\XXX_{(\Delta,\Delta^0)}$ the toric stack associated to $(\Delta,\Delta^0)$.
Let $A^*(\XXX_{(\Delta,\Delta^0)})$ denote the
integral Chow ring of $\XXX_{(\Delta,\Delta^0)}$.
Suppose that rays in $\Delta$ span the vector space $N\otimes_{\ZZ}\RR$.
Then there exists an isomorphism of graded rings
\[
\textup{SR}_{(\Delta,\Delta^0)}\stackrel{\sim}{\to}A^*(\XXX_{(\Delta,\Delta^0)})\]
which sends $\DDD_{\rho}$ to
$c_1(\OO_{\XXX_{(\Delta,\Delta^0)}}(\DDD_{\rho}))$.
If cones $\sigma, \tau\in\Delta$ span $\gamma\in \Delta$,
then we have
\begin{equation}
\tag*{(*)}
  [\VVV(\sigma)]\cdot [\VVV(\tau)]=[\VVV(\gamma)],
\end{equation}
in $A^*(\XXX_{(\Delta,\Delta^0)})$.
\end{Theorem}

As the first step to Theorem~\ref{main2},
we will calculate the Picard goup of a toric stack.

\begin{Proposition}
\label{pic}
Let $(\Delta,\Delta^0)$ be a stacky fan and
$\XXX_{(\Delta,\Delta^0)}$ the toric stack associated to $(\Delta,\Delta^0)$.
Let $\Pic (\XXX_{(\Delta,\Delta^0)})$ denote the group
of invertible sheaves on $\XXX_{(\Delta,\Delta^0),\textup{\'et}}$.
Suppose that rays in $\Delta$ span the vector space $N\otimes_{\ZZ}\RR$.
Then there exists a natural isomorphism of groups
\[
\bigoplus_{\rho\in \Delta(1)}\ZZ\cdot \DDD_{\rho}/(\Sigma_{\rho\in\Delta(1)}\langle m,n_{\rho}\rangle\cdot \DDD_{\rho})_{m\in M}\stackrel{\sim}{\longrightarrow}\Pic (\XXX_{(\Delta,\Delta^0)}),
\]
\[
\DDD_{\rho}\mapsto \OO_{\XXX_{(\Delta,\Delta^0)}}(\DDD_{\rho})
\]
where $\bigoplus_{\rho\in \Delta(1)}\ZZ\cdot \DDD_{\rho}$
is a free abelian group generated by $\{ \DDD_{\rho}\}_{\rho\in \Delta(1)}$
\textup{(}we abuse notation\textup{)}.
\end{Proposition}

\Proof
Put $G:=G_{(\Delta,\Delta^0)}$.
Observe first that every invertible sheaf on $L_{\Delta}$
is trivial.
Indeed $L_{\Delta}$ is a smooth toric variety,
and any invertible sheaf (line bundle) $\LL$ on $L_{\Delta}$
is represented by a linear form of torus-invariant
divisors with integer coefficients.
Every torus-invariant divisor on $L_{\Delta}$ comes from
some toric divisor on $\AAAA^{\Delta(1)}$
(recall $L_{\Delta}\subset \AAAA^{\Delta(1)}$ !).
Since every torus-invariant divisor on $\AAAA^{\Delta(1)}$
is a principal divisor, thus $\LL$ is trivial.
Then we obtain
\begin{Claim}
The Picard group $\Pic (\XXX_{(\Delta,\Delta^0)})$
on $\XXX_{(\Delta,\Delta^0)}$
%$\cong \FF_{(\Delta,\Delta^0)}\cong[L_{\Delta}/G_{(\Delta,\Delta^0)}]$
is isomorphic to
\[
\bigoplus_{\rho\in \Delta(1)}\ZZ\cdot \DDD_{\rho}/(\Sigma_{\rho\in\Delta(1)}\langle m,n_{\rho}\rangle\cdot \DDD_{\rho})_{m\in M}.
\]
\end{Claim}

{\it Proof of Claim.}
Let $\LL$ be an invertible sheaf on $\XXX_{(\Delta,\Delta^0)}$.
Let $p:L_{(\Delta,\Delta^0)}\to 
\XXX_{(\Delta,\Delta^0)}\cong \FF_{(\Delta,\Delta^0)}\cong [L_{\Delta}/G_{(\Delta,\Delta^0)}]$ be a projection and set $M:=p^*\LL$.
By the above observation, $M$ is a trivial invertible sheaf.
Let $\textup{pr}_1:L_{\Delta}\times G\to L_{\Delta}$
denote the first projection and $a:L_{\Delta}\times G\to L_{\Delta}$
the action.
By the descent theory, the invertible sheaf $\LL$ on $\XXX_{(\Delta,\Delta^0)}\cong[L_{\Delta}/G]$ amounts exactly to
an isomorphism of invertible sheaves $\sigma:\textup{pr}_1^*M\stackrel{\sim}{\to}a^*M$ which satisfies the usual cocycle condition.
Notice that the complement $\AAAA^{\Delta(1)}-L_{\Delta}$ has codimension
more than 1. Thus $\Aut (\textup{pr}_1^*M)$ is canonically isomorphic to
\[
\Hom_{(k\textup{-schemes})}
(L_{\Delta}\times G,\GG_m)
\cong \Hom_{(k\textup{-schemes})}(\AAAA^{\Delta(1)}\times G,\GG_m)
\cong \Hom_{(k\textup{-schemes})}(G,\GG_m)
\]
where $\Aut (\textup{pr}_1^*M)$ is the group
of automorphisms of the invertible sheaf $\textup{pr}_1^*M$.
Let $(M,\sigma_0:\textup{pr}_1^*M\stackrel{\sim}{\to}a^*M)$
be a pair corresponding to a trivial invertible sheaf on $\XXX_{(\Delta,\Delta^0)}$.
Taking account of the cocycle condition,
there exists an isomorphism of groups
\[
\Hom_{(\textup{group}\ k\textup{-schemes})}(G,\GG_m)\to \Pic (\XXX_{(\Delta,\Delta^0)})
\]
which sends $u\in \Hom_{(\textup{group}\ k\textup{-schemes})}(G,\GG_m)\subset
H^0(G,\OO^*_{G})$
to $(M, \textup{pr}_1^*M\stackrel{u}{\to}\textup{pr}_1^*M\stackrel{\sigma_0}{\to}a^*M)$. Here $u:\textup{pr}_1^*M\to\textup{pr}_1^*M$
is defined by $1\mapsto u$.
On the other hand, there exists an isomorphism of groups
\[
\Hom_{(\textup{group $k$-schemes})}(G,\GG_{m})\cong \bigoplus_{\rho\in \Delta(1)}\ZZ\cdot \DDD_{\rho}/(\Sigma_{\rho\in\Delta(1)}\langle m,n_{\rho}\rangle\cdot \DDD_{\rho})_{m\in M}.
\]
Hence we obtain our claim.
\QED

Let us go back to the proof of Proposition.
If $C$ denotes the complement $\XXX_{(\Delta,\Delta^0)}-\Spec k[M]$
there exists the excision sequence (\cite[Proposition 2.4.1]{KR})
\[
A_{d-1}(C)\to A_{d-1}(\XXX_{(\Delta,\Delta^0)})\to A_{d-1}(\Spec k[M])\to 0
\]
where $A_{d-1}(\bullet)$ denote integral Chow groups (cf. \cite{KR},\cite{EG}).
Moreover $A_{d-1}(\Spec k[M])=0$ and
$A_{d-1}(C)=\ZZ^{\oplus (\textup{irreducible components of }C)}=\bigoplus_{\rho\in \Delta(1)}\ZZ\cdot \DDD_{\rho}$.
Hence we have a surjective map $\bigoplus_{\rho\in \Delta(1)}\ZZ\cdot \DDD_{\rho}\to A_{d-1}(\XXX_{(\Delta,\Delta^0)})\cong \Pic (\XXX_{(\Delta,\Delta^0)})$ (the last isomorphism follows
from \cite[Proposition 18]{EG}).
Observe that for any $m\in M$, $\Sigma_{\rho\in\Delta(1)}\langle m,n_{\rho}\rangle\cdot \DDD_{\rho}$ maps to zero.
To see this, consider a free abelian group $\tilde{N}:=\bigoplus_{\rho\in \Delta(1)}\ZZ\cdot e_{\rho}$,
and define a homomorphism of abelian group $h:\tilde{N}\to N$
by $e_{\rho}\mapsto n_{\rho}$.
Let $\tilde{\Delta}\subset \tilde{N}\otimes_{\ZZ}\RR$
be the sub-fan of $\oplus_{\rho\in \Delta(1)}\RR_{\ge0}\cdot e_{\rho}$
such that a cone $\gamma\in\oplus_{\rho\in \Delta(1)}\RR_{\ge0}\cdot e_{\rho}$
lies in $\tilde{\Delta}$ if and only if $h_{\RR}(\gamma)$ is a cone in $\Delta$.
Then by \cite[Corollary 3.9]{I2}
the associated morphism of stacky fans
$h:(\tilde{\Delta},\tilde{\Delta}_{\textup{can}}^{0})\to (\Delta,\Delta^0)$
determines a smooth surjective morphism
$p:X_{\tilde{\Delta}}\to \XXX_{(\Delta,\Delta^0)}$.
Note that the morphism $h:\tilde{\Delta}\to \Delta$
of fans corresponds to the composite $X_{\tilde{\Delta}}\stackrel{p}{\to} \XXX_{(\Delta,\Delta^0)}\stackrel{\pi_{(\Delta,\Delta^0)}}{\longrightarrow} X_{\Delta}$.
By the construction, a rational function $m\in M$ of $\XXX_{(\Delta,\Delta^0)}$
defines a torus-invariant divisor $p^*(m)=\Sigma_{\rho}\langle m,n_{\rho}\rangle\cdot E_{\rho}$ on $X_{\tilde{\Delta}}$, where $E_{\rho}$ is the reduced torus-invariant divisor corresponding to $\RR_{\ge0}\cdot e_{\rho}\in \tilde{\Delta}$.
Since $p^{-1}(\DDD_{\rho})=E_{\rho}$, thus a rational function
$m\in M$ induces $\Sigma_{\rho}\langle m,n_{\rho}\rangle\cdot \DDD_{\rho}$
on $\XXX_{(\Delta,\Delta^0)}$.
%To see this, notice that we can regard $m$ as a rational function on $\Spec k[M%]$
%and it defines the divisor
%$\Sigma_{\rho}\langle m,v_{\rho}\rangle\cdot D_{\rho}$
%on $X_{\Delta}$, where $v_{\rho}$ is the first lattice point of $\rho$.
%If $\pi_{(\Delta,\Delta^0)}:\XXX_{(\Delta,\Delta^0)}\to X_{\Delta}$
%denotes the coarse moduli map,
%then $\pi_{(\Delta,\Delta^0)}^{-1}(D_{\rho})=l_{\rho}\cdot \DDD_{\rho}$
%where $l_{\rho}$ is the level of $\Delta^0$ on $\rho$.
%Therefore $m\in M$ defines a divisor
%$\Sigma_{\rho}l_{\rho}\cdot \langle m,v_{\rho}\rangle\cdot \DDD_{\rho}
%=\Sigma_{\rho}\langle m,n_{\rho}\rangle\cdot \DDD_{\rho}$.
Therefore we have a surjective map
\[
\bigoplus_{\rho\in \Delta(1)}\ZZ\cdot \DDD_{\rho}/(\Sigma_{\rho\in\Delta(1)}\langle m,n_{\rho}\rangle\cdot \DDD_{\rho})_{m\in M}\longrightarrow\Pic (\XXX_{(\Delta,\Delta^0)}),\ \DDD_{\rho}\mapsto c_1(\OO(\DDD_{\rho})).
\]
To complete the proof it suffices to show that this map is injective.
The group $\Pic (\XXX_{(\Delta,\Delta^0)})$
is isomorphic to $\bigoplus_{\rho\in \Delta(1)}\ZZ\cdot \DDD_{\rho}/(\Sigma_{\rho\in\Delta(1)}\langle m,n_{\rho}\rangle\cdot \DDD_{\rho})_{m\in M}$.
Thus the injectiveness follows since
$\bigoplus_{\rho\in \Delta(1)}\ZZ\cdot \DDD_{\rho}/(\Sigma_{\rho\in\Delta(1)}\langle m,n_{\rho}\rangle\cdot \DDD_{\rho})_{m\in M}$ is a finitely generated abelian group.
\QED

Now we will prove Theorem~\ref{main2}.

{\it Proof of Theorem~\ref{main2}.}
First of all, if two cones $\sigma$ and $\tau$ in $\Sigma$
span the cone $\gamma$,
then $\VVV(\sigma)$ and $\VVV(\tau)$ intersect transversally
at $\VVV(\gamma)$ by \cite[Proposition 4.19]{I}
(it also follows from the quotient presentation).
This implies the relation (*).
Next we shall show that the map from Stanley-Reisner ring
to the Chow ring is an isomorphism.
By Corollary~\ref{quot},
the toric stack $\XXX_{(\Delta,\Delta^0)}$ has the quotient presentation $[L_{\Delta}/G_{(\Delta,\Delta^0)}]$,
where $G_{(\Delta,\Delta^0)}$ is isomorphic to $\GG_m^k\times H$.
Here $H$ is a finite abelian group, and $k=\#\Delta(1)-\dim (\XXX_{(\Delta,\Delta^0)})$.
Set $G=G_{(\Delta,\Delta^0)}$ and $W:=L_{(\Delta,\Delta^0)}$.
Note that 
$W$ has the form $\mathbb{A}^{\Delta(1)}-Z$ where $Z$ is a close subscheme.
Decompose $H$ into $\ZZ / m_1\ZZ \times \cdots \times\ZZ / m_r\ZZ$.
Note that we may view $G$ as the closed subgroup of the maximal 
algebraic torus in $W$.
According to the definition of Chow groups due to Edidin-Graham (\cite{EG}),
first of all, in order to compute the Chow ring of $\XXX_{(\Delta,\Delta^0)}$,
we shall construct a certain
$N(k+r)$-dimensional representation of $G$, i.e.,
an action of $G$ on the affine space $V:=\mathbb{A}^{N(k+r)}$
such that $V$ has an open set $U$ on which
$G$ acts freely and whose complement has 
codimension more than $N-1$.
To this aim,
by choosing the primitive $m_i$-th root $\zeta_i$
in the base field for $1\le i\le r$,
we embed $G=\GG_m^k\times \ZZ / m_1\ZZ \times \cdots\times \ZZ / m_r\ZZ$
into the closed group subscheme
of $\GG_m^k\times\GG_m^r$ as follows:
\[
G=\GG_m^k\times \ZZ / m_1\ZZ \times \cdots \times\ZZ / m_r\ZZ
\hookrightarrow \GG_m^k\times\GG_m^r,
\]
\[
(u,l_1,\ldots,l_r)\mapsto (u,\zeta_1^{l_1},\ldots,\zeta_r^{l_r}),
\]
where $u \in \GG_m^k$ and $(l_1,\ldots,l_r)\in \ZZ / m_1\ZZ \times \cdots \times\ZZ / m_r\ZZ$.
This embedding yields the action of $G$ on an affine space
$\mathbb{A}^{k+r}$.
We extend this action to $\mathbb{A}^{k+r}\times\cdots\times \mathbb{A}^{k+r}=\mathbb{A}^{N(k+r)}$ diagonally.
Set $U_i:=\mathbb{A}^{(i-1)(k+r)}\times\GG_m^{(k+r)}\times\mathbb{A}^{(N-i)(k+r)}$, and $U:=\cup_{1\le i\le N}U_i$.
Then the action of $G$ on $U$ is free and $\mathbb{A}^{N(k+r)}-U$ has codimension more than $N-1$.
Then we
have $A_i(\XXX_{(\Delta,\Delta^0)})=A_{i+N(k+r)}((W\times U)/G)$
for $N>\dim (\XXX_{(\Delta,\Delta^0)})-i$ (cf. \cite{EG}, \cite{KR}).
Here $G$ acts on $W\times U$ diagonally (this is a free action).

Next, let us compute the Chow group of $(W\times U)/G$.
To this aim, we show that $(W\times U)/G$ is {\it a smooth toric variety}.
First, notice that since $G$ is a linear reductive algebraic group,
and thus
by geometric invariant theory, we see that
$(W\times U)/G$ is an algebraic variety.
Moreover $W\times U$ and $G$ are smooth,
$(W\times U)/G$ is smooth.
$W\times U$ contains the maximal torus $\GG_m^{\Delta(1)}\times\GG_m^{N(k+r)}$
as a dense open subset, and this torus naturally acts on $W\times U$.
Since
$\GG_m^{\Delta(1)}\times\GG_m^{N(k+r)}$
is commutative, the action of $\GG_m^{\Delta(1)}\times\GG_m^{N(k+r)}$
on $W\times U$ descends to $(W\times U)/G$.
Thus the action of $(\GG_m^{\Delta(1)}\times\GG_m^{N(k+r)})/G$
on itself is naturally extended to $(W\times U)/G$.
Recall the fact that quotients of diagonalizable groups are diagonalizable.
Since $(\GG_m^{\Delta(1)}\times\GG_m^{N(k+r)})/G$ is connected,
$(\GG_m^{\Delta(1)}\times\GG_m^{N(k+r)})/G$
is an algebraic split torus. Hence we conclude that
$(W\times U)/G$ is a smooth toric variety.

Let $N'=\ZZ^{d+N(k+r)}$ be a lattice and $\Sigma$ the fan in $N'\otimes_{\ZZ}\RR=N'_{\RR}$
such that the associated toric variety $X_{\Sigma}$
is $(W\times U)/G$.
Let us denote by $\Sigma(1)'$ the set of torus-invariant divisors on
$(W\times U)/G$
arising from the torus-invariant divisors on $W$ (or equivalently $[W/G]$).
(Notice that $\Sigma(1)'$ also bijectively corresponds to the set of torus-invariant
divisors on the toric variety $\mathbb{A}^{\Delta(1)}$.)
The projection $(W\times U)/G \to [W/G]$ is a torus-equivariant
morphism and induces a bijective map 
\[
q:\Delta(1) \to \Sigma(1)',
\]
via the flat pull-back of $(W\times U)/G\to [W/G]$.
The following is a key Lemma for the proof.

\begin{Lemma}
\label{tech}

\renewcommand{\theenumi}{(\roman{enumi})}
\begin{enumerate}

\item The fan $\Sigma$ is a non-singular fan,
and the rays in $\Sigma$ span the vector space $N_{\RR}'$.

\item Let  $\alpha_1,\ldots,\alpha_a$ be rays in $\Sigma(1)'$,
and $ \beta_1,\ldots ,\beta_b$ be rays in $\Sigma(1)-\Sigma(1)'$.
Suppose that the cones $\langle \alpha_1,\ldots,\alpha_a \rangle$
and $\langle  \beta_1,\ldots ,\beta_b \rangle$ are in $\Sigma$.
Then the cone $\langle \alpha_1,\ldots,\alpha_a,\beta_1,\ldots,\beta_b \rangle$
is in  $\Sigma$.

\item If $\beta_1,\ldots,\beta_{N-1}$ is rays in $\Sigma(1)-\Sigma(1)'$,
then $ \beta_1,\ldots,\beta_{N-1} $ span
a cone in $\Sigma$.

\end{enumerate}
\end{Lemma}

\Proof
(i). Since $(W\times U)/G$ is smooth, thus $\Sigma$ is a
non-singular fan.
To show that the rays in $\Sigma$ span the vector space $N_{\RR}'$,
notice that $W\times U$ is also a toric variety and
the projection $W\times U \to (W\times U)/G$ is a
torus-equivariant map.
We denote by $\tilde{\Sigma}$ the fan in $N_{\RR}''$ 
which represents the toric variety $W\times U$.
Let us denote by $\phi:(\tilde{\Sigma},N'')\to (\Sigma,N')$
the map of fans which represents the projection
$W\times U \to (W\times U)/G$.
Since the projection is surjective,
$\phi$ induces the surjective map $\phi_{\RR}:N_{\RR}''\to N_{\RR}'$.
On the other hand, the rays in $\tilde{\Sigma}$ span the
vector space $N_{\RR}''$ and thus our claim follows.

(ii). By the construction of action of $G$ on $W\times U$,
the torus-invariant divisors on $W\times U$ bijectively correspond
to the torus-invariant divisors on $(W\times U)/G$.
Thus we have $\tilde{\Sigma}(1)\cong \Sigma(1)$.
Let
$\tilde{\alpha_1},\ldots,\tilde{\alpha_a}$ ,
and $ \tilde{\beta_1},\ldots ,\tilde{\beta_b}$ be rays in $\tilde{\Sigma}$
which correspond to $\alpha_1,\ldots,\alpha_a$,
and $ \beta_1,\ldots ,\beta_b$ respectively.
Both of sets of rays
$\{ \tilde{\alpha_1},\ldots,\tilde{\alpha_a} \}$,
and $\{ \tilde{\beta_1},\ldots ,\tilde{\beta_b} \}$
span cones in $\tilde{\Sigma}$,
and the cycles $V(\langle \tilde{\alpha_1},\ldots,\tilde{\alpha_a} \rangle)$
and $V(\langle \tilde{\beta_1},\ldots,\tilde{\beta_b} \rangle)$
intersect.
Therefore $V(\langle {\alpha_1},\ldots,{\alpha_a} \rangle)$
and $V(\langle {\beta_1},\ldots,{\beta_b} \rangle)$
intersect and thus 
 ${\alpha_1},\ldots,{\alpha_a}, {\beta_1},\ldots,{\beta_b}$
span the cone in $\Sigma$.

(iii). Let $\tilde{\beta}_1,\ldots,\tilde{\beta}_{N-1}$ be
rays in $\tilde{\Sigma}$ which correspond to
$\beta_1,\ldots,\beta_{N-1}$.
These rays $\tilde{\beta}_1,\ldots,\tilde{\beta}_{N-1}$
also correspond to the torus-invariant divisors
in $U$.
Since $U=\cup_{1\le i\le N}(\mathbb{A}^{(i-1)(k+r)}\times\GG_m^{(k+r)}\times\mathbb{A}^{(N-i)(k+r)}) \subset \mathbb{A}^{N(k+r)}$,
$V(\langle\tilde{\beta}_1,\ldots,\tilde{\beta}_{N-1}\rangle)$
is a non-empty set.
Thue the rays $\tilde{\beta}_1,\ldots,\tilde{\beta}_{N-1}$
span the cone $\langle \tilde{\beta}_1,\ldots,\tilde{\beta}_{N-1} \rangle$
in $\tilde{\Sigma}$.
Since $\tilde{\Sigma}(1)\cong \Sigma(1)$ and $\Sigma$ is non-singular,
$\phi_{\RR}(\langle \tilde{\beta}_1,\ldots,\tilde{\beta}_{N-1} \rangle)
=\langle {\beta}_1,\ldots,{\beta}_{N-1} \rangle$
and this completes the proof.
\QED

From now on the positive integer $N$ is assumed to be greater than 
$1$.
The group $A_{d-1+N(k+r)}([(W\times U)/G])=A_{d-1+N(k+r)}(X_{\Sigma})$
is generated by torus-invariant divisors $V(\rho)$
as $\rho$ ranges over $\Sigma(1)$ (cf. \cite[Proposition 2.1]{FS}).
The rational equivalent relations on these divisors are
generated by the linear forms $\Sigma_{\rho \in \Sigma(1)}\langle m,v_{\rho}\rangle\cdot V(\rho)$ as $m$ ranges over $M'=\Hom (N',\ZZ)$.
Here $v_{\rho}$ is the first lattice point of $\rho$.
We denote by $I$ the subgroup of $\oplus_{\rho \in \Sigma(1)}\ZZ\cdot D_{\rho}$
generated by the above linear forms $\Sigma_{\rho \in \Sigma(1)}\langle m,v_{\rho}\rangle\cdot D_{\rho}$.
On the other hand, by Proposition~\ref{pic}, 
there exists a natural isomorphism
\[
A_{d-1+N(k+r)}([(W\times U)/G])
=A_{d-1}(\XXX_{(\Delta,\Delta^0)}) \cong
\bigoplus_{\rho\in \Delta(1)}\ZZ\cdot\DDD_{\rho}/(\Sigma_{\rho\in\Delta(1)}\langle m,n_{\rho}\rangle\cdot \DDD_{\rho},
m \in M).
\]
Thus $A_{d-1+N(k+r)}([(W\times U)/G])$ is generated by
the torus-invariant divisors arising from the torus-invariant
divisors on $\XXX_{(\Delta,\Delta^0)}\cong [W/G]$.
Therefore we can choose the element 
\begin{equation}
\tag*{(*)}
D_{\xi}-\Sigma_{\rho \in \Sigma(1)'}a_{\rho}\cdot D_{\rho}\ \ \ (a_{\rho} \in \ZZ)
\end{equation}
in $I$ for each ray $\xi$ in $(\Sigma(1)-\Sigma(1)')$.
The Chow ring $A^{*}([(W\times U)/G])$
is generated by torus-invariant divisors $D_{\rho}$
as $\rho$ ranges over $\Sigma(1)'$.
Indeed, notice that $[(W\times U)/G]$ is smooth. Thus
by \cite[Proposition 2.1]{FS},
 the Chow ring $A^{*}([(W\times U)/G])$
is generated by torus-invariant divisors $D_{\rho}$
as $\rho$ ranges over $\Sigma(1)'$
because $A_{d-1+N(k+r)}([(W\times U)/G])$ is generated by
the torus-invariant divisors arising from the torus-invariant
divisors on $\XXX_{(\Delta,\Delta^0)}\cong [W/G]$.
Let $\rho_1,\ldots,\rho_k$ be rays in $\Delta$.
Then $\langle \rho_1,\ldots,\rho_k \rangle $ is a cone in $\Delta$
if and only if $\langle q(\rho_1),\ldots,q(\rho_k) \rangle $ 
is a cone in $\Sigma$.
The element 
$\Sigma_{\rho\in\Delta(1)}\langle m,n_{\rho}\rangle\cdot D_{q(\rho)}$
lies in $I_{\Sigma}$.
Therefore there exist surjective homomorphisms of graded rings
\[
\xi:\ZZ[\DDD_{\rho},\rho\in \Delta(1)]/(I_{(\Delta,\Delta^0)}+J_{\Delta})
\to \ZZ[D_{\rho},\rho\in \Sigma(1)]/(I_{\Sigma}+J_{\Sigma}),
\]
where $\xi(\DDD_{\rho})=D_{q(\rho)}$, and
\[
\eta:\ZZ[D_{\rho},\rho\in \Sigma(1)]/(I_{\Sigma}+J_{\Sigma})
\to A^*([(W\times U)/G]).
\]
where $\eta(D_{\rho})=V(\rho)$.
%Here $\ZZ[D_{\rho},\rho\in \Delta(1)]/(I_{\Delta}+J_{\Delta})$
%is the Stanley-Reisner ring of $\Delta\subset N'_{\RR}$.

\vspace{2mm}

Next we show the following Proposition.

\begin{Proposition}
The composition $\eta\circ\xi$
is bijective modulo degree $N$.
\end{Proposition}
\Proof
To prove this, we first show that
$\xi$ is bijective modulo degree $N$.
Let $J_{\Sigma}^{(1)}$ (resp. $J_{\Sigma}^{(2)}$) be the ideal
generated by the monomials $D_{\rho_1}\cdots D_{\rho_v}$
such that  $\rho_1,\ldots,\rho_v$ are in $\Sigma(1)'$ 
(resp. $\Sigma(1)-\Sigma(1)'$)
and $\langle \rho_1,\ldots,\rho_v \rangle \notin \Sigma(1)$.
Then by Lemma~\ref{tech} (ii), $J_{\Sigma}$ is generated
by two ideals $J_{\Sigma}^{(1)}$ and $J_{\Sigma}^{(2)}$. 
Moreover by Lemma~\ref{tech} (iii),
$J_{\Sigma}^{(2)}$ is generated by the monomials 
whose degree are greater than $N-1$.
On the other hand, by Proposition~\ref{pic} and \cite[Proposition 2.1]{FS},
$\xi$ is bijective modulo degree $2$.
Since $I_{\Sigma}$ are generated by linear forms
and contains elements (*), we see that $\xi$
is bijective modulo degree $N$.
Finally, Proposition follows from the next lemma.

\begin{Lemma}
\label{fulst}
The natural morphism of graded rings
\[
\eta:\ZZ[D_{\rho}, \rho\in \Sigma(1)]/(I_{\Sigma}+J_{\Sigma})
{\to}
A^*([(W\times U)/G])=A^*(X_{\Sigma}),
\]
is an isomorphism.
\end{Lemma}

\Proof
Let $\tau$ be an $s$-dimensional cone in $\Sigma$.
Let $m$ be an element in $\tau^{\perp}\cap M'$.
We define $e(\tau,m)$ to be 
\begin{center}
$\displaystyle \sum_{\langle \rho_1,\ldots,\rho_{s+1}\rangle \supset\tau \atop
\operatorname{in} \Sigma} \langle m,n_{\langle \rho_1,\ldots,\rho_{s+1}\rangle,\tau}\rangle D_{\rho_1}\cdots D_{\rho_{s+1}}$.
\end{center}
Here $n_{\langle \rho_1,\ldots,\rho_{s+1}\rangle,\tau}$ is
a lattice point in $\langle \rho_1,\ldots,\rho_{s+1}\rangle$
whose image generates the one dimensional lattice
$N_{\langle \rho_1,\ldots,\rho_{s+1}\rangle}'/N_{\tau}'$.
Here $N'_{\langle \rho_1,\ldots,\rho_{s+1}\rangle}$
(resp. $N'_{\tau}$) is a sublattice of $N'$
generated by $\langle \rho_1,\ldots,\rho_{s+1}\rangle\cap N'$
(resp. $\tau\cap N'$).
Let us denote by
$v_{\langle \rho_1,\ldots,\rho_{s+1}\rangle,\tau}$ the first
lattice point of the ray of $\langle \rho_1,\ldots,\rho_{s+1}\rangle$,
 which is not contained in $\tau$ (it is unique since $\Sigma$ is a non-singular fan, in particular simplicial).
 Since $\Sigma$
 is a {\it non-singular} fan, $e(\tau,m)$ is equal to
 \begin{center}
$\displaystyle \sum_{\langle \rho_1,\ldots,\rho_{s+1}\rangle \supset\tau \atop
\operatorname{in} \Sigma} \langle m,v_{\langle \rho_1,\ldots,\rho_{s+1}\rangle,\tau}\rangle D_{\rho_1}\cdots D_{\rho_{s+1}}$.
\end{center}

 To show the Lemma, firstly,
 note that
 the ring $\eta:\ZZ[D_{\rho}, \rho\in \Sigma(1)]/(I_{\Sigma}+J_{\Sigma})$
 has no element whose degree is greater than $\dim [(W\times U)/G]$.
 Indeed, by Lemma~\ref{tech} (i), it follows from an elementary argument
 as in the proof of \cite[10.7.1]{D}.
 Then taking account into \cite[Proposition 2.1 (b)]{FS},
 what we have to show is the following claim.

\begin{Claim}
\label{co}
Let $(I_{\Sigma}+J_{\Sigma})_k$ be the $\ZZ$-module
generated by the homogeneous elements of the ideal $I_{\Sigma}+J_{\Sigma}$
whose degree are equal to $k$. Suppose that $k\le \dim [(W\times U)/G]$.
Let $(L_{\Sigma})_k$ be the $\ZZ$-module generated by elements $e(\tau,m)$
$(\tau \in \Sigma(k-1), m\in \tau^{\perp}\cap M')$
and
the monomials $D_{\rho_1}\cdots
D_{\rho_k}$ such that $\langle \rho_1,\ldots,\rho_k \rangle \notin \Sigma$.
Then there exists the canonical homomorphism
\begin{align*}
\bigl( \displaystyle \bigoplus_{ \rho_1,\ldots,\rho_k  \in \Sigma(1) \atop \rho_i \neq \rho_j \operatorname{if} i\neq j} \ZZ D_{\rho_1}\cdots
D_{\rho_k}\bigr) \bigl/(L_{\Sigma})_k &\to
\bigl( \displaystyle \bigoplus_{\rho_1,\ldots,\rho_s \in \Sigma(1) \atop r_1,\ldots,r_s \in \ZZ_{\ge 1},\ r_1+\cdots+r_s=k }\ZZ D_{\rho_1}^{r_1}\cdots D_{\rho_s}^{r_s}\bigl)\bigl/(I_{\Sigma}+J_{\Sigma})_k, \\
\end{align*}
and it is an isomorphism.
\end{Claim}

\Proof 
Firstly, we shall show that $(I_{\Sigma}+J_{\Sigma})_k $ contains
$(L_{\Sigma})_k$.
To this aim, put $\tau=\langle \rho_1,\ldots,\rho_{k-1}\rangle$.
Then for $m$ in $\tau^{\perp}\cap M'$, we have
\begin{align*}
e(\tau,m) &\equiv D_{\rho_1}\cdots D_{\rho_{k-1}}\bigl(\displaystyle\sum_{\rho\in \Sigma(1)-\{
\rho_1,\ldots,\rho_{k-1} \} \atop \operatorname{such\ that} \ \langle \rho,\tau \rangle \in \Sigma}\langle m,v_{\rho} \rangle D_{\rho}\bigl) \\
&\equiv D_{\rho_1}\cdots D_{\rho_{k-1}} \bigl(\displaystyle\sum_{\rho\in \Sigma(1)}\langle m,v_{\rho} \rangle D_{\rho}\bigl).
\end{align*}
Here by $\equiv$, we mean ``modulo monomials $D_{\alpha_1}\cdots
D_{\alpha_k}$ such that $\langle \rho_1,\ldots,\rho_k \rangle \notin \Sigma$".
Clearly, if $\langle \alpha_1,\ldots,\alpha_k\rangle \notin \Sigma$,
$D_{\alpha_1}\cdots D_{\alpha_k}$ is contained in $(I_{\Sigma}+J_{\Sigma})_k$.
Hence  $(I_{\Sigma}+J_{\Sigma})_k $ contains
$(L_{\Sigma})_k$.
Taking account into Lemma~\ref{tech} (i),
by an elementary argument as in the proof of \cite[Lemma 10.7.1]{D},
we see that
the canonical injective map
\begin{align*}
( \displaystyle \bigoplus_{ \rho_1,\ldots,\rho_k  \in \Sigma(1) \atop \rho_i \neq \rho_j \operatorname{if} i\neq j} \ZZ D_{\rho_1}\cdots
D_{\rho_k}\bigr) \bigl/ \bigl((I_{\Sigma}+J_{\Sigma})_k\cap ( \displaystyle \bigoplus_{ \rho_1,\ldots,\rho_k  \in \Sigma(1) \atop \rho_i \neq \rho_j \operatorname{if} i\neq j} \ZZ D_{\rho_1}\cdots
D_{\rho_k}\bigr)\bigr) \\
\to 
\bigl( \displaystyle \bigoplus_{\rho_1,\ldots,\rho_s \in \Sigma(1) \atop r_1,\ldots,r_s \in \ZZ_{\ge 1},\ r_1+\cdots+r_s=k }\ZZ D_{\rho_1}^{r_1}\cdots D_{\rho_s}^{r_s}\bigl)\bigl/(I_{\Sigma}+J_{\Sigma})_k \\
\end{align*}
is bijective.
To prove our claim, it suffices only to prove
that $(L_{\Sigma})_k$ contains
\[
(I_{\Sigma}+J_{\Sigma})_k^{\circ}:=(I_{\Sigma}+J_{\Sigma})_k \cap ( \displaystyle \bigoplus_{ \rho_1,\ldots,\rho_k  \in \Sigma(1) \atop \rho_i \neq \rho_j \operatorname{if} i\neq j} \ZZ D_{\rho_1}\cdots
D_{\rho_k}\bigr).
\]
Notice that $(L_{\Sigma})_k$ generates
rational equivalent relations on
the $(\dim [(W\times U)/G]-k)$-dimensional torus-invariant cycles (cf. \cite[Proposition 2.1 (b)]{FS}).
Then our assertion is clear since $(I_{\Sigma}+J_{\Sigma})_1$ defines rational equivalent relations
and if $\langle \rho_1,\ldots,\rho_k \rangle \notin \Sigma$,
the intersection $V(\rho_1)\cap \cdots\cap V(\rho_k)$ is empty.
Thus we have $(I_{\Sigma}+J_{\Sigma})_k^{\circ} \subset (L_{\Sigma})_k$.
This implies our claim.
\QED

Now let us go back to the proof of Theorem~\ref{main2} again.

Note that $A_i(\XXX_{(\Delta,\Delta^0)})=A_{i+N(k+r)}((W\times U)/G)$
for $N>\dim \XXX_{(\Delta,\Delta^0)}-i$ (cf. \cite{EG}).
The map
\[
\eta\circ \xi:\ZZ[\DDD_{\rho},\rho\in \Delta(1)]/(I_{(\Delta,\Delta^0)}+J_{\Delta})
\to A^*([(W\times U)/G])
\]
is bijective modulo degree $N$.
Since $N$ is an arbitrary positive integer,
thus we complete the proof of Theorem~\ref{main2}.
\QED

\begin{Example}
\label{exam}
\renewcommand{\theenumi}{(\roman{enumi})}

Let us look at some examples of toric stacks.

\begin{enumerate}
\begin{comment}
\item
Let $N\cong \ZZ\cdot e_1\oplus \ZZ\cdot e_2$ be a lattice.
Let $\Sigma$ be a complete fan in $N_{\RR}$
with rays $\rho_1=\RR_{\ge 0}\cdot(e_1+e_2)$,
$\rho_2=\RR_{\ge 0}\cdot(-e_1+e_2)$, and $\rho_3=\RR_{\ge 0}\cdot(-e_1-2e_2)$.
Let us denote by $A^*(\XXX_{\Sigma})$ the integral Chow ring
of the associated toric stack $\XXX_{\Sigma}$.
Then by Theorem~\ref{Main2}, there exists the isomorphism of graded rings
\[
\ZZ[D_1,D_2,D_3]/(D_1-D_2-D_3, D_1+D_2-2D_3, D_1D_2D_3)\stackrel{\sim}{\to}A^*(\XXX_{\Sigma})
\]
where $D_i$ maps to $\VVV(\rho_i)$ in $A^1(\XXX_{\Sigma})$ for $i=1,2,3$.
Thus there exist isomorphisms of groups
$A^k(\XXX_{\Sigma})\cong \ZZ$ for $k=0,1,2$
and $A^k(\XXX_{\Sigma})\cong \ZZ/6\ZZ$ for $k\ge 3$.
\end{comment}

\item Let $\Delta$ be a complete fan in $N_{\RR}$
with rays $\rho_1=\RR_{\ge 0}\cdot(2e_1-e_2)$,
$\rho_2=\RR_{\ge 0}\cdot(-e_1+2e_2)$, and $\rho_3=\RR_{\ge 0}\cdot(-e_1-e_2)$.
Then by Theorem~\ref{main2} there exists an isomorphism of graded rings
\begin{center}
$\textup{SR}_{\Delta}=\ZZ[D_1,D_2,D_3]/(2D_1-D_2-D_3, -D_1+2D_2-D_3, D_1D_2D_3)\stackrel{\sim}{\to}A^*(\XXX_{(\Delta,\Delta_{\textup{can}}^0)})$,
\end{center}
where $D_i$ maps to $c_1(\OO(\VVV(\rho_i)))$ in $A^1(\XXX_{(\Delta,\Delta_{\textup{can}}^0)})$ for $i=1,2,3$.
Thus we see that $A^0(\XXX_{(\Delta,\Delta_{\textup{can}}^0)})\cong \ZZ$, $A^1(\XXX_{(\Delta,\Delta_{\textup{can}}^0)})\cong \ZZ\oplus \ZZ/3\ZZ$,
 $A^2(\XXX_{(\Delta,\Delta_{\textup{can}}^0)})\cong \ZZ\oplus (\ZZ/3\ZZ)^{\oplus 2}$,
 and $A^k(\XXX_{(\Delta,\Delta_{\textup{can}}^0)})\cong (\ZZ/3\ZZ)^{\oplus 3}$ for $k\ge 3$.
 On the other hand, the Chow group (resp. operational Chow group)
 $A_*(X_{\Delta})$ (resp. $A^*(X_{\Delta})$) of the
 toric variety $X_{\Delta}$ are computed as follows:
$ A_0(X_{\Delta})=\ZZ,\ A_1(X_{\Delta})=\ZZ\oplus\ZZ/3\ZZ,\ 
 A_2(X_{\Delta})=\ZZ,\ A_k(X_{\Delta})=0\ \operatorname{for}\ k\ge 3
 \operatorname{or} k<0$,
and
$A^k(X_{\Delta})=\ZZ \ \operatorname{for}\ 0\le k\le 2 ,\ 
A^k(X_{\Delta})=0 \ \operatorname{for}\ k\ge 3$,
by using \cite[Proposition 2.1 and Proposition 2.4]{FS}.

\item
(Weighted projective line).
Let $N=\ZZ\cdot e$ be a lattice of rank one.
Let $(\Sigma,\Sigma^0)$ be a complete fan in $N_{\RR}$
with rays $\rho_1=\RR_{\ge 0}\cdot e$
and $\rho_2=\RR_{\ge 0}\cdot (-e)$.
Let $a\cdot e$ and $b\cdot (-e)$
be generators of $\Sigma^0$ on $\rho_1$ and $\rho_2$ respectively,
where $a$ and $b$ are positive integers.
If $a$ and $b$ are coprime, 
we have an isomorphisms of stacks
$\PPP(a,b)=[(\AAAA^2-(0,0))/\GG_m]\cong \XXX_{(\Sigma,\Sigma^0)}$
(the proof is left to the readers).
Here the action of $\GG_m$ is defined by $u\cdot(x,y)\mapsto (u^ax,u^by)$
for $u\in\GG_m$.
Then we have an isomorphism
\[
\ZZ[t]/(ab\cdot t^2)\cong A^*(\XXX_{(\Sigma,\Sigma^0)})\cong A^*(\PPP(a,b)).
\]

\end{enumerate}
\end{Example}


\begin{thebibliography}{99}


%\bibitem{EGA1}
%J. Dioudonn\'{e} and A. Grothendieck,
%\'Elements de G\'{e}om\'{e}trie Alg\'{e}briques I,
%Publ. Math. IHES 4 (1960)


%\bibitem{EGA3}
%J. Dioudonn\'{e} and A. Grothendieck,
%\'Elements de G\'{e}om\'{e}trie Alg\'{e}briques Chapitre III,
%Publ. Math. IHES 11 (1961)


%\bibitem{AV}
%D. Abramovich and A. Vistoli,
%Compactifying the space of stable maps,
%Jour. Amer. Math. Soc. 15 (2002) no.1 p27--75.

%\bibitem{AMRT}
%A. Ash, D. Mumford, M. Rapoport, and Y.-S. Tai,
%Smooth compactifications of locally symmetric varieties, Math. Sci. Press,
%Bookline, MA, (1975)



\bibitem{C}
D. Cox,
The functor of a smooth toric variety,
Tohoku Math. J. (2) 47 (1995) no.2. 251--262.




\bibitem{C2}
D. Cox,
The homogeneous coordinate ring of a toric variety,
J. Alg. Geom. 4 (1995) 17--50.

\bibitem{D}
V. I. Danilov,
The geometry of toric varieties,
Uspekhi Mat. Nauk. 33:2 (1978) 85--134.


\bibitem{EG}
D. Edidin and W. Graham,
Equivariant intersection theory,
Invent. Math. 131, 595--634 (1998).




%\bibitem{FI}
%W. Fulton,
%Intersection Theory,
%Springer-Verlag, (1984).


\bibitem{F}
W. Fulton,
Introduction to Toric Varieties,
Ann. Math. Stud. Princeton Univ. Press (1993).

\bibitem{FS}
W. Fulton and B. Sturmfels,
Intersection theory of toric varieties,
Topology 36 (1997) 335--353.


%\bibitem{G}
%H. Gillet,
%Intersection theory on algebraic stacks and Q-varieties,
%J. Pure Appl. Algebra 34 (1984) 193--240.


\bibitem{I}
I. Iwanari,
Logarithmic geometry, minimal free resolutions and toric algebraic stacks,
preprint.





\bibitem{I2}
I. Iwanari,
The category of toric stacks,
preprint (2006),
Arxiv: math.AG/0610548

\bibitem{KM}
S. Keel and S. Mori,
Quotients by groupoids,
Ann. Math. 145 (1997), 193--213.



\bibitem{KR}
A. Kresch,
Cycle groups for Artin stacks,
Invent. Math. 138 (1999) 495--536




\bibitem{LM}
G. Laumon and L. Moret-Bailly,
Champs Alg\'ebriques,
Springer-Verlag (2000).


\bibitem{M}
D. Mumford, J. Fogarty and F. Kirwan,
Geometric Invariant Theory,
Third edition Springer-Verlag (1992).






\bibitem{Hom}
M. Olsson,
Hom-stacks and restriction of scalars,
Duke Math. J. 134 (2006) no.1 139--164.


\bibitem{TO}
B. Totaro,
The Chow ring of a classfying space,
in Algebraic K-Theory, Proceedings of Symposia in Pure Math. v. 67, Amer. Math. Soc. (1999), 249-281. 


\bibitem{V}
A. Vistoli,
Intersection theory on algebraic stacks and their moduli spaces,
Invent. Math. 97 (1989), 613--670.

\bibitem{V2}
A. Vistoli,
The Chow ring of $\mathcal{M}_2$,
Invent. Math. 131 (1998) 635--644.




\end{thebibliography}
\end{document}